\documentclass[11pt,a4paper]{amsart}

\usepackage{amsfonts, amsmath, amssymb, amsthm}
\theoremstyle{plain}
\newtheorem{thm}{Theorem}[section]
\newtheorem{lem}[thm]{Lemma}
\newtheorem{prop}[thm]{Proposition}
\newtheorem{cor}[thm]{Corollary}

\theoremstyle{definition}
\newtheorem{defn}[thm]{Definition}
\newtheorem{exmp}[thm]{Example}

\newtheorem*{claim}{Claim}

\theoremstyle{remark}
\newtheorem{rem}[thm]{Remark}

\usepackage{txfonts}
\usepackage{mathbbol}
\usepackage{ae}

\usepackage{color}
\usepackage{graphicx}
\usepackage{url}
\usepackage[abs]{overpic}
\usepackage{subfigure}
\usepackage{booktabs}

\newcommand\cA{{\mathcal A}}

\newcommand\cT{{\mathcal T}}

\newcommand\RR{{\mathbb R}}
\newcommand\ZZ{{\mathbb Z}}
\newcommand\abs[1]{|#1|}

\newcommand\SetOf[2]{\left\{#1\,\vphantom{#2}\right.\left|\vphantom{#1}\,#2\right\}}
\newcommand\smallSetOf[2]{\{#1\,|\,#2\}}
\newcommand\conv{\operatorname{conv}}
\newcommand\Sym{\operatorname{Sym}}
\newcommand\bbc[1]{\operatorname{bbc}({#1})}
\newcommand\ncp[2]{\operatorname{ncp}_{#1}\!\left({#2}\right)}
\newcommand\mirror[1]{\operatorname{M}(#1)}
\newcommand\fissure[3]{\operatorname{fis}_{#1}(#2,#3)}
\newcommand\cube[1]{\operatorname{C}_{#1}}
\newcommand\simplex[1]{\Delta_{#1}}
\newcommand\cyclic[2]{\operatorname{cyc}_{#1}({#2})}

\renewcommand{\phi}{\varphi}

\begin{document}

\title{Neighborly Cubical Polytopes and Spheres}
\author[Joswig and Schr\"oder]{Michael Joswig \and Thilo Schr\"oder}
\address{Michael Joswig, Fachbereich Mathematik, AG~7, TU Darmstadt, 64289 Darmstadt, Germany}
\email{joswig@mathematik.tu-darmstadt.de}
\address{Thilo Schr\"oder, Institut f\"ur Mathematik, MA 6-2, TU Berlin, 10623 Berlin, Germany}
\email{thilosch@math.tu-berlin.de}
\thanks{Both authors are supported by Deutsche Forschungsgemeinschaft, DFG Research Group ``Polyhedral Surfaces.''}
\date{\today}

\begin{abstract}
  We prove that the neighborly cubical polytopes studied by G\"unter M. Ziegler and the first author~\cite{MR1758054}
  arise as a special case of the neighborly cubical spheres constructed by Babson, Billera, and Chan~\cite{MR1489110}.
  By relating the two constructions we obtain an explicit description of a non-polytopal neighborly cubical sphere and,
  further, a new proof of the fact that the cubical equivelar surfaces of McMullen, Schulz, and Wills~\cite{MR657865}
  can be embedded into~$\RR^3$.
\end{abstract}

\maketitle

\section{Introduction}

Our point of departure is a paper by Babson, Billera, and Chan~\cite{MR1489110}, in which the authors introduce an
inductive construction of cubical $d$-spheres from certain sequences of simplicial $(d-1)$-balls and their boundary
spheres of dimension $d-2$.  It turns out that such cubical spheres reflect many properties of the simplicial spheres
involved, but in a cubical disguise.  In particular, this way the boundary of a neighborly simplicial $(d-1)$-polytope
yields a neighborly cubical $d$-sphere, that is, a sphere which has the same $\lfloor (d-1)/2\rfloor$-skeleton as some
high-dimensional cube.  Later it was shown by entirely different methods that there even exist neighborly cubical
spheres which are polytopal~\cite{MR1758054}.

Our first result, Theorem~\ref{thm:GeneralizedBBC}, establishes a non-recursive combinatorial description of the cubical
spheres studied by Babson, Billera, and Chan.  From this it can be inferred that particular sequences of pulling
triangulations of cyclic polytopes yield the polytopal spheres studied in~\cite{MR1758054}, see
Corollary~\ref{cor:NCSeqNCP}.  As a further benefit of this direct description we observe in Theorem~\ref{thm:polytopal}
that the simplicial spheres involved in the construction necessarily must be polytopal in order to yield polytopal
cubical spheres.  In Theorem~\ref{thm:Altshuler} we derive that the neighborly cubical $5$-sphere on $2048$ vertices
that we explicitly construct from a certain $3$-sphere $\mathcal{M}^{10}_{425}$, found by Altshuler~\cite{MR0467531} and
whose non-polytopality was proved by Bokowski and Garms~\cite{MR919873}, cannot be polytopal.

A seemingly different topic is the construction of polyhedral surfaces of `unusually large genus,' pioneered by
Coxeter~\cite{C:SkewPolyhedra} and Ringel~\cite{MR0075586} and continued by McMullen, Schulz, and
Wills~\cite{MR657865,MR727027}.  Yet we can show that the neighborly cubical $4$-polytopes of~\cite{MR1758054} contain
cubical surfaces with $n$ vertices of genus $\mathcal{O}(n\log n)$ in their boundary.  Via Schlegel diagrams this gives
a simple new proof of the known fact that there are such surfaces which can be embedded into~$\RR^3$ with straight
faces.  This answers a question of G\"unter M. Ziegler~\cite{ZieglerPC}.

\section{Cubical and Simplicial Complexes}

We review the ingredients of the construction of Babson, Billera, and Chan~\cite{MR1489110} with slight generalizations.

The reader is advised to consult the monographs \cite{MR1976856} and \cite{MR1311028} for general information about
convex polytopes and related topics.

\subsection{Cubes}
Consider the $2d$ affine halfspaces
\[ H_i^+=\SetOf{x\in\RR^d}{x_i\le 1}\ \text{and}\ H_i^-=\SetOf{x\in\RR^d}{-x_i\le 1}, \]
which define the $d$-dimensional cube
\[ \cube{d}=H_1^+\cap\dots\cap H_d^+\cap H_1^-\cap\dots\cap H_d^-. \]
The intersections
\[ F_i^\sigma=\partial H_i^\sigma\cap\cube{d}=\SetOf{x\in\cube{d}}{\sigma x_i=1},\qquad \sigma \in \{\pm 1 \}, \]
of their boundaries with the cube are precisely the \emph{facets} of $\cube{d}$, that is, its maximal proper faces.  The
\emph{vertices} of $\cube{d}$ are all the $2^d$ vectors of length~$d$ with coordinates $\pm 1$, and their joint convex
hull is the cube $\cube{d}$.  Two facets $F_i^+$ and $F_i^-$ are parallel, and hence they do not share any vertex of
$\cube{d}$.

It is a general fact about convex polytopes that each \emph{proper face}, that is, any intersection of the polytope with
a supporting hyperplane, can be written as the intersection of facets.  Thus, for the special case of the $d$-cube, each
non-empty face $F$ can be written uniquely as $F=F_{i_1}^{\sigma_{i_1}}\cap\dots F_{i_k}^{\sigma_{i_k}}$, where
$i_j\ne i_l$ for any $j\ne l$.  By letting $\sigma_j=0$ for all $j\not\in\{i_1,\dots,i_k\}$ the non-empty face $F$ can
be identified with the ordered sequence $(\sigma_1,\dots,\sigma_d)$ of signs $+1,0,-1$ of length~$d$.  Conversely,
each such sign vector defines a face. For ease of notation we often omit the $1$'s of $\pm1$.

The intersection of a $k$-dimensional cube face, or $\emph{$k$-face}$ for short, with a facet is either empty or a
$(k-1)$-face.  This readily implies that the dimension of the face $(\sigma_1,\dots,\sigma_d)$ equals the number of
$0$-entries in its sign vector.

\begin{figure}[htbp]\centering
  \includegraphics[width=.75\textwidth]{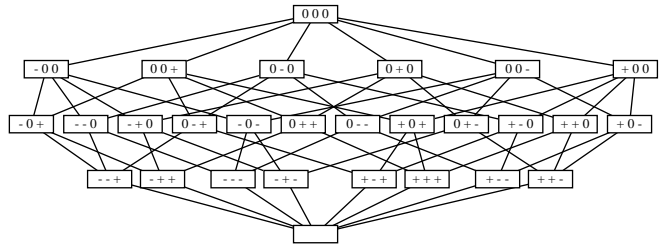}
  \caption{Face lattice of $\cube{3}$.}
\end{figure}

The reflections at the coordinate hyperplanes $\smallSetOf{x\in\RR^d}{x_i=0}$ generate an elementary abelian group
$\Sigma_d$ of order $2^d$ which acts sharply transitively on the vertices of $\cube{d}$.  The full automorphism group
$\Gamma_d$ of the $d$-cube is isomorphic to a semi-direct product of $\Sigma_d$ and the stabilizer of a vertex, which is
$\Sym_d$, the symmetric group of degree~$d$.  In fact, $\Gamma_d$ is the wreath product $\ZZ_2\wr\Sym_d$.

\subsection{Regular Cell Complexes and Posets}
\label{sec:PosetsAndCellComplexes}

A \emph{regular cell complex} is a family $\mathcal{C}$ of closed balls in a Hausdorff space $X_{\mathcal{C}}$, such
that the interiors of the balls partition $X_{\mathcal{C}}$ and the boundary of each ball is the union of balls in
$\mathcal{C}$. The elements of $\mathcal{C}$ are called cells.  The topology of a regular cell complex is completely
determined by its face poset; see Bj\"orner~\cite[Prop.~3.1]{MR746039}.

An \emph{abstract $d$-complex} $\mathcal{P}$ is a ranked poset of rank $d+1$ such that there exist unique lower bounds
for any set of elements (that is, $\mathcal{P}$ is a meet semi-lattice) and every order ideal is combinatorially
isomorphic to the face lattice of a polytope.  The elements of this partially ordered set are called faces.  The
\emph{boundary} of a finite abstract $d$-complex is the set of faces of rank $d$ contained in only one maximal face. An
abstract $d$-complex is a \emph{cubical} (\emph{simplicial}) \emph{complex} if every face of~$\mathcal{P}$ is
combinatorially isomorphic to a cube (simplex).  These are CW posets of polyhedral type as introduced by
Bj\"orner~\cite{MR746039}. Thus we are able to speak of topological properties of abstract $d$-complexes. An abstract
$2$-complex representing a connected $2$-manifold without boundary is a \emph{polyhedral surface}.

\subsection{Mirroring}
An abstract simplicial complex $\Delta$ on $d$ vertices can be seen as a subcomplex of the $(d-1)$-dimensional simplex.
The cube $\cube{d}$ is a \emph{simple} $d$-polytope, that is, each of its vertex figures is a $(d-1)$-simplex.  Here the
\emph{vertex figure} of a polytope $P$ at a vertex~$v$ is the intersection of $P$ with an affine hyperplane which
separates $v$ from the other vertices of~$P$.  We construct a (cubical) subcomplex of $\cube{d}$ which corresponds to a
simultaneous embedding of $\Delta$ into the vertex figures of all the vertices of $\cube{d}$ such that these embeddings
are invariant under the action of the group $\Sigma_d$.

Following Babson, Billera, and Chan~\cite{MR1489110} we encode $\Delta$ in a non-standard way: If $1,2,\dots,d$ are the
vertices of~$\Delta$ we associate with a face $\phi\in\Delta$ the characteristic function of its complement
$\{1,\dots,d\}\setminus\phi$.  Using this description we define
\[ 
\mirror{\Delta}= 
\SetOf{(\sigma_1,\dots,\sigma_d) \in \{0,\pm1\}^d}{(\abs{\sigma_1},\dots,\abs{\sigma_d})\in\Delta},
\]
the \emph{mirror complex} of $\Delta$, which is a subcomplex of the $d$-cube.

\begin{figure}[htbp]\centering
  \label{fig:mirrorcomplex}
  \raisebox{2em}{
    \begin{overpic}[width=.25\textwidth]{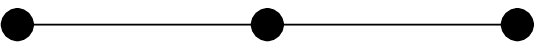}
      \put(0,9){$1$}
      \put(42,9){$2$}
      \put(85,9){$3$}
    \end{overpic}
  }\qquad
  \includegraphics[width=.3\textwidth]{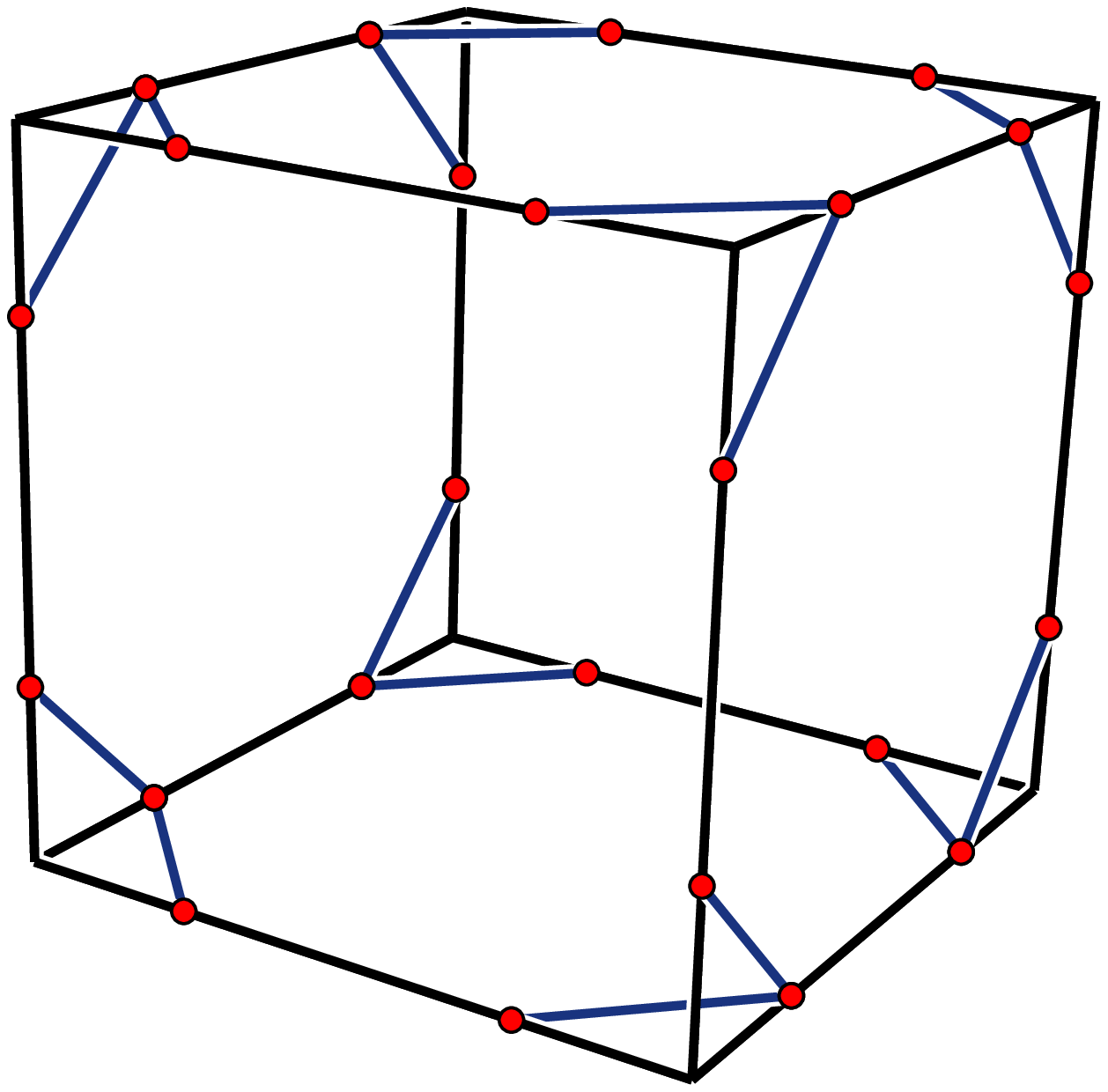}\
  \includegraphics[width=.3\textwidth]{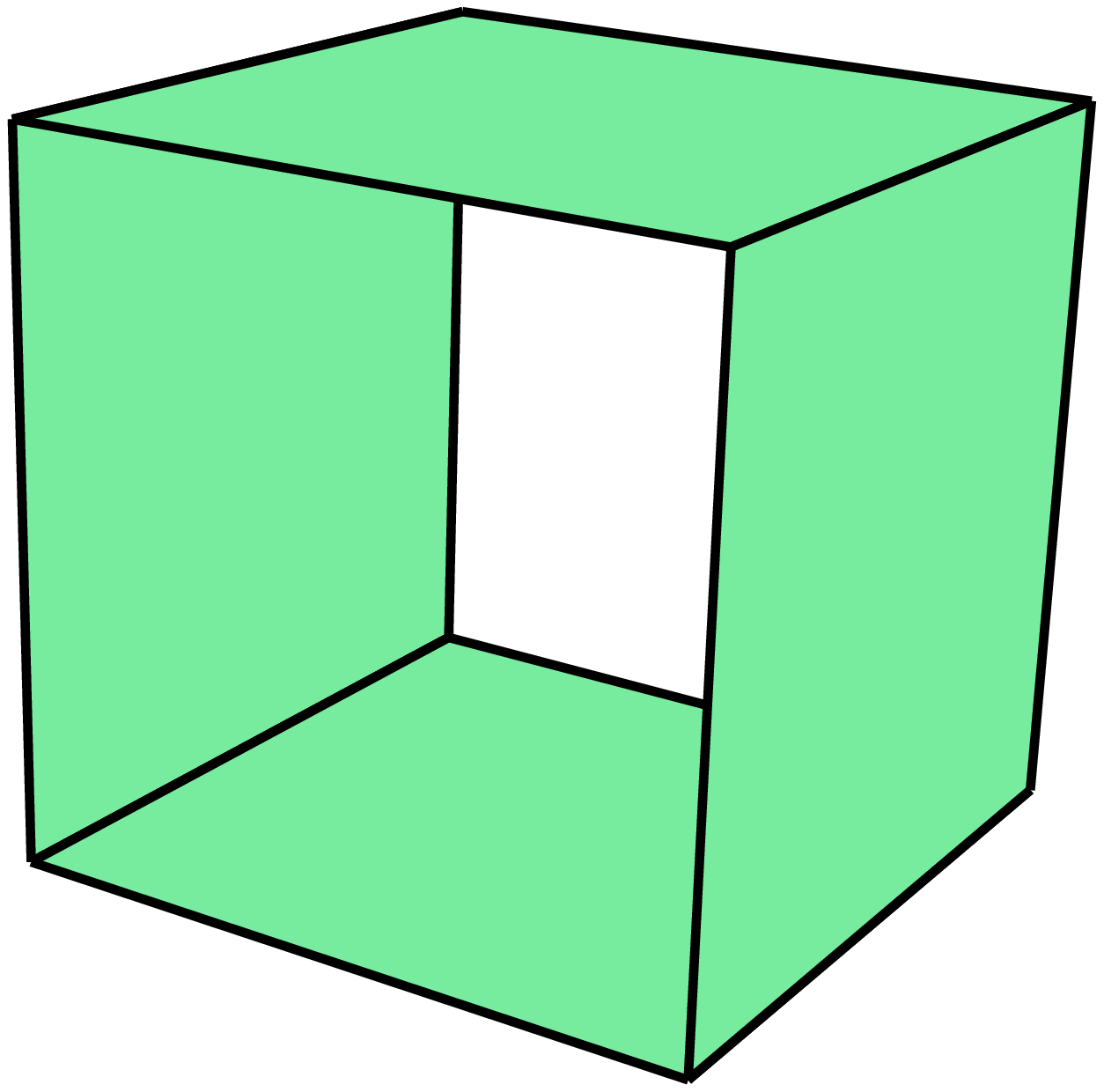}

  \caption{Simplicial complex $\Delta$ on three vertices, its symmetric embedding into $\cube{d}$, and the mirror
    complex $\mirror{\Delta}$.}
\end{figure}

By construction the vertex figure of the mirror complex of a simplicial complex is isomorphic to the simplicial complex
itself. Thus if the simplicial complex is a $d$-sphere, then its mirror complex is a $(d+1)$-manifold. The $f$-vector of
the mirror complex of a simplicial complex $\Delta$ on $n$ vertices is $f_i(\mirror{\Delta}) = 2^{n-i} f_{i-1}(\Delta)$.
Further the mirror complex of boundary of the simplicial complex is the boundary of the mirror complex:
$\mirror{\partial \Delta} = \partial \mirror{\Delta}$.

\begin{prop}
  Let $\Delta$ be a simplicial complex with automorphism group~$\Gamma$.  Then the automorphism group of
  $\mirror{\Delta}$ is isomorphic to the semi-direct product $\Sigma_d\rtimes\Gamma$.
\end{prop}

\subsection{Fissuring}

\label{sec:Fissuring}
The cubical fissure or fissuring is an operation that produces a new cubical complex from a given one. Let $C$ be a pure
cubical $d$-complex, $C_1$ and $C_2$ facet-disjoint $d$-dimensional subcomplexes of $C$ such that $C_1 \cup C_2 = C$.
The \emph{cubical fissure} $\fissure{C}{C_1}{C_2}$ of $C$ between $C_1$ and $C_2$ is defined by lifting $C_1$ to height
one, dropping $C_2$ to height minus one and filling in the fissure with $C_1 \cap C_2 \times [-1,1]$.  The corresponding
poset is
\[
\fissure{C}{C_1}{C_2} = (C_1 \times \{+1\}) \cup (C_1 \cap C_2 \times \{0\}) \cup (C_2 \times \{-1\}).
\]
with $+1 < 0$ and $-1 < 0$ in the last component. If $C$ is a subcomplex of the $n$-cube given as sign vectors of length
$n$, then the cubical fissure canonically yields a subcomplex of the $(n+1)$-cube. The cubical fissure between a
subcomplex $C_1$ and its complement $C_2$ yields a complex consisting of the subcomplex $C_1$ and its complement
connected via a prism over the boundary of $C_1$.

\begin{exmp}
  Recall the simplicial complex $\Delta$ of Figure~\ref{fig:mirrorcomplex} on three vertices. The mirror complex $C_1 =
  \mirror{\Delta}$ is a subcomplex of the boundary complex of the $3$-cube $C = \partial{\cube{3}}$. So the cubical
  fissure between $C_1$ and its complement in $C$ is a subcomplex of the boundary of the $4$-cube.
  \begin{figure}[htbp]\centering
  \label{fig:fissuring}
  \includegraphics[width=.3\textwidth]{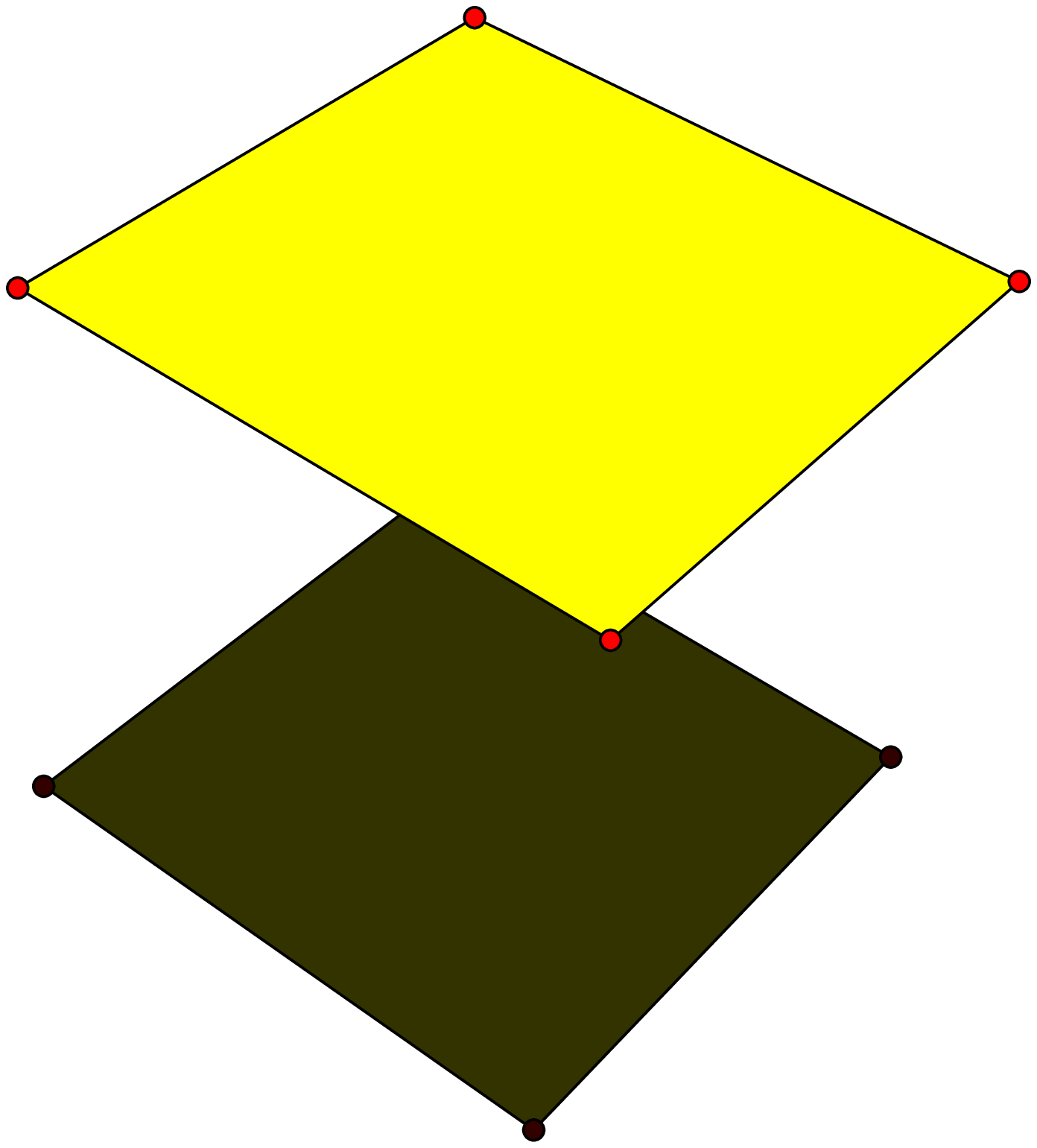}\
  \includegraphics[width=.3\textwidth]{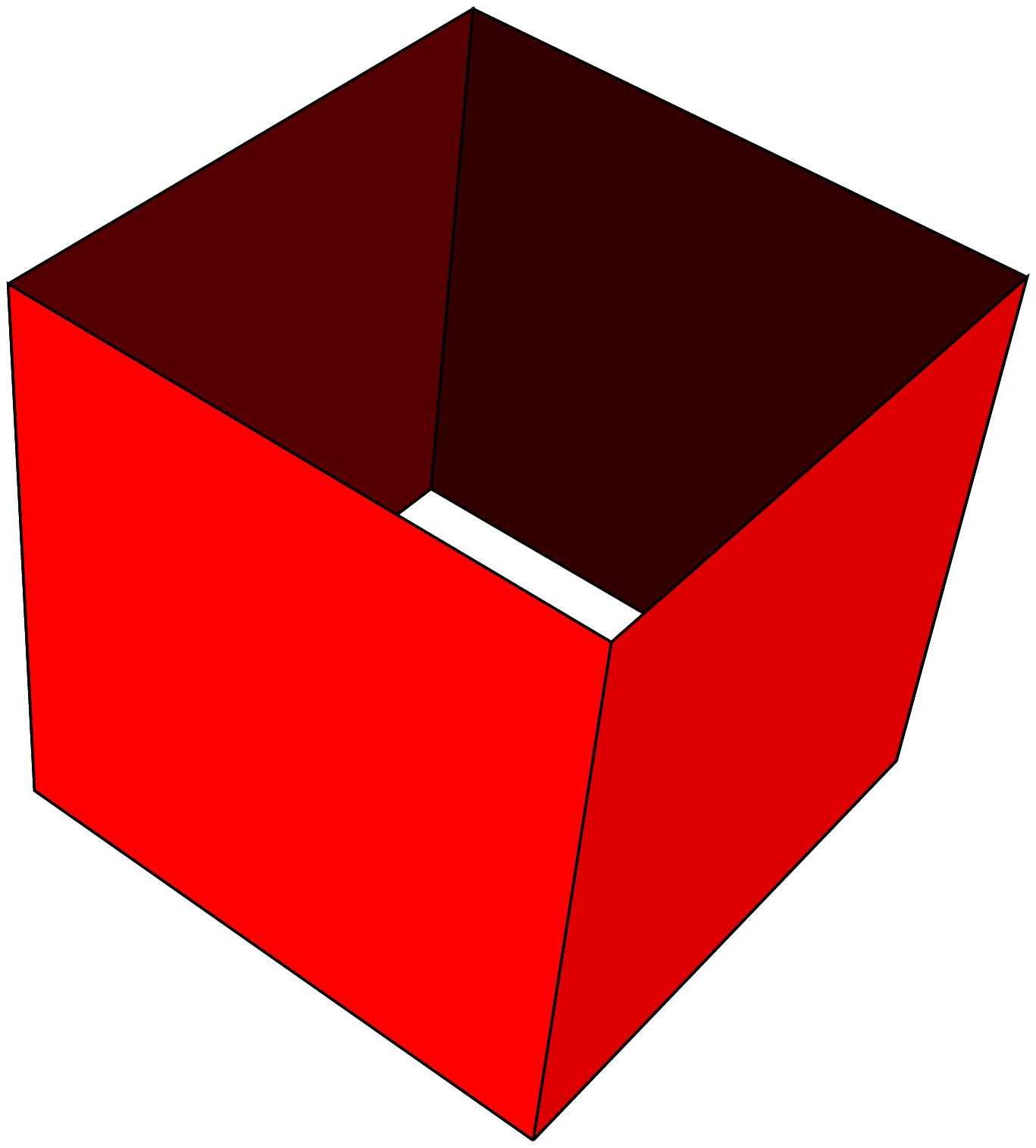}\
  \includegraphics[width=.3\textwidth]{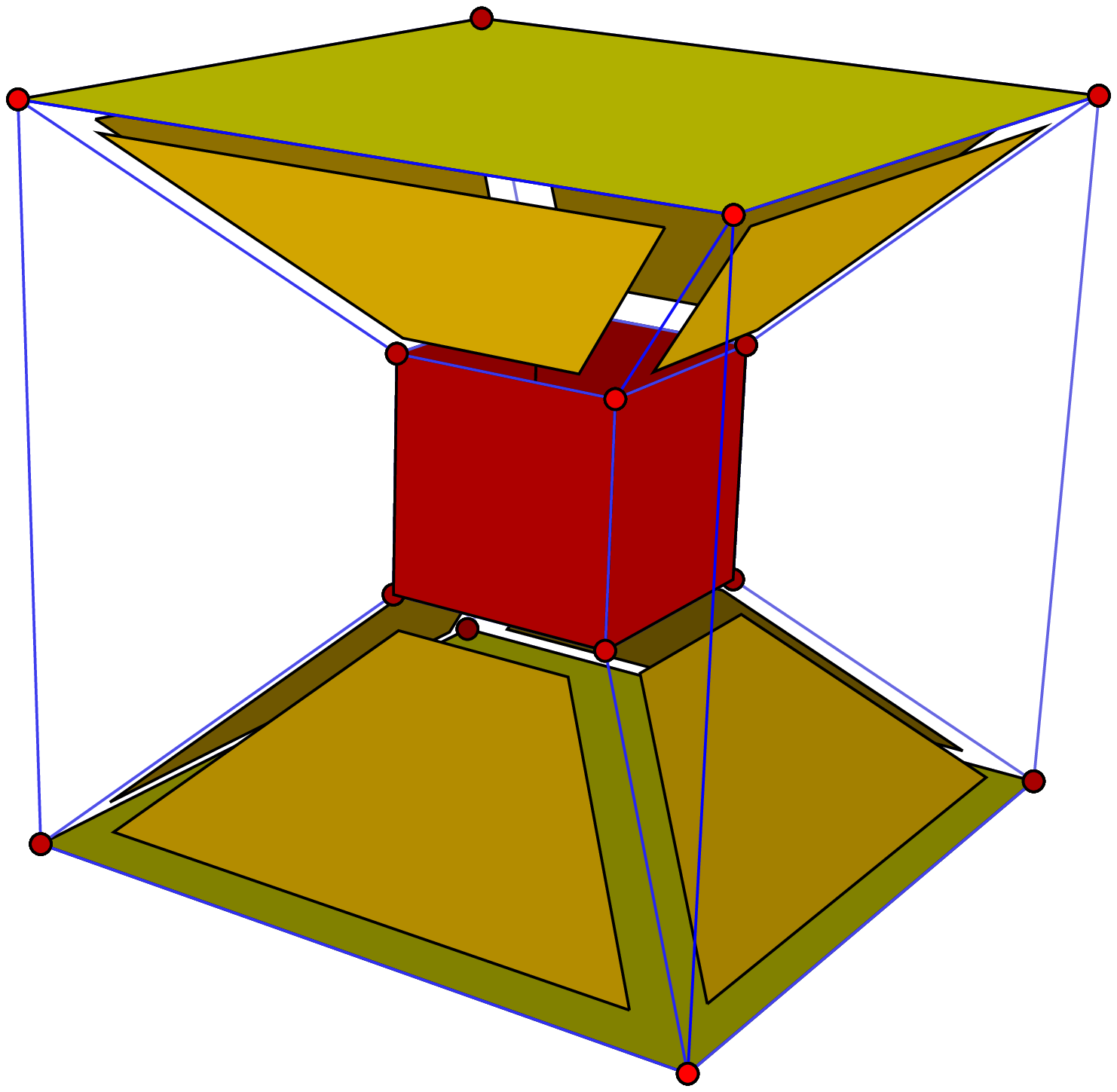}
  \caption{Two subcomplexes of the $3$-cube and their cubical fissure in the Schlegel diagram of the $4$-cube.}
\end{figure}
\end{exmp}

\subsection{BBC Sequences of Simplicial Balls}
Babson, Billera, and Chan~\cite{MR1489110} proved the existence of neighborly cubical spheres.  Their approach is based
on an inductive construction using triangulations of cyclic polytopes, mirror complexes and fissuring.  A close
inspection of their proof motivates the following definition.  While this does look a bit technical it is one of
the keys to our main results.

\begin{defn}
\label{defn:BBCSequences}
Let $V = \{v_0,\ldots,v_{n}\}$ be a set of vertices, $B_i$ be a simplicial $(d-1)$-ball on the vertex set
$V_i=\{v_0,\ldots,v_{i}\}$ for $i=d,\ldots,n-1$ with $d \geq 2$.  Then $(T_i)_{i=d+1}^{n}$ is a \emph{BBC sequence} of
simplicial $d$-balls if
  \[
  T_i = B_{i-1} * v_{i} \quad\text{and}\quad B_{i-1} \subseteq \partial T_{i-1}\quad\text{for $i=d+1,\ldots,n$}.
  \]
\end{defn}

\begin{rem}
  \label{rem:BBCSequences}
  Let $P$ be a simplicial $(d+1)$-polytope with vertices $v_1,\ldots,v_n$. We may assume that the vertices of $P$ lie in
  general position since $P$ is simplicial. So every subset of at least $d+2$ of its vertices is again the vertex set of
  a simplicial $(d+1)$-polytope. Thus by taking an arbitrary ordering $(v_{l_1}, \ldots, v_{l_n})$ of the vertices of
  $P$ we obtain a sequence of $(d+1)$-balls $P_i = \conv({v_{l_1}, \ldots, v_{l_i}})$ for $i = d+2,\ldots,n$. The
  corresponding BBC sequence $(T_i)_{i=d+2}^n$ is the sequence of pulling triangulations of $P_i$ with respect to
  $v_{l_i}$.
\end{rem}

Note that $(T_i)_{i=d+1}^{n}$ is a BBC sequence if and only if each boundary sphere $\partial T_i$ is \emph{directly
  obtainable} from its predecessor $\partial T_{i-1}$ in the sense of Altshuler~\cite{MR0397744}.  In loc. cit. he
proved that there is a not-directly obtainable $3$-sphere on ten vertices.  Remark~\ref{rem:BBCSequences} implies that
such a sphere cannot be polytopal.

\begin{rem}
  Since each sphere of dimension at most~$2$ is polytopal it follows that every BBC sequence is a sequence of pulling
  triangulations of polytopes.  But there exist BBC sequences of simplicial $4$-balls such that the final boundary
  sphere is not polytopal.  One such sequence is described in Section~\ref{sec:NonPolytopalSphere},
  Table~\ref{tab:BBCSequenceOfAlthulerSphere}.
\end{rem}

Let $\cT=(T_i)_{i=d}^{n-1}$ be a BBC sequence of simplicial $(d-1)$-balls.  We inductively define cubical complexes
$S_k$ for $k=d+1,\ldots,n$ with $2^k$ vertices.  Babson, Billera, and Chan begin their inductive definition with $S_d$,
which is two $d$-cubes identified at their complete boundary. But since this does not yield a regular cell complex,
we start with $S_{d+1}=\partial \cube{d+1}$ the boundary of the $(d+1)$-cube (that is, the mirror complex of the
boundary of the $d$-simplex). Then for $k = d+2,\ldots,n$ we recursively define:
\begin{align*}
  S_{k} & = \fissure{S_{k-1}}{\mirror{T_{k-1}}}{S_{k-1}\setminus\mirror{T_{k-1}}} \\
  &= \left(\mirror{T_{k-1}} \times \{+1\}\right) \cup \left(\partial\mirror{T_{k-1}} \times \{0\}\right) \cup \left((S_{k-1}
    \setminus \mirror{T_{k-1}}) \times \{-1\}\right)
\end{align*}
We denote the final cubical complex $S_n$ by $\bbc{\cT}$.

\begin{thm}
  \label{thm:BBC}
  Let $\cT=(T_i)_{i=d}^{n-1}$ be a BBC sequence of simplicial $(d-1)$-balls.  Then the cubical complexes $S_k$ and, in
  particular, $\bbc{\cT}$ are cubical $d$-spheres.
\end{thm}
\begin{proof}
  This is part of what is proved by Babson, Billera, Chan~\cite[Theorem~3.1]{MR1489110}.  The essential step is to see
  that the mirror complex of $T_k$ is actually a subcomplex of $S_k$, such that fissuring is possible.
\end{proof}

If the simplicial balls $T_i$ are equipped with a piecewise-linear (PL) structure, then the resulting cubical spheres
are also PL~\cite{MR1489110}.

Babson, Billera, and Chan were interested in neighborly cubical spheres.  The corresponding construction is included
in the following section.

\subsection{Neighborliness}

A simplicial complex is \emph{(simplicially) $k$-neighborly} if every $k$-subset of its vertices is a face. In other
words, its $(k-1)$-skeleton is isomorphic to the $(k-1)$-skeleton of a simplex. A simplicial $d$-sphere is
\emph{(simplicially) neighborly} if it is simplicially $\lfloor(d+1)/2\rfloor$-neighborly. The same is defined in the
cubical case. A cubical complex is \emph{(cubically) $k$-neighborly} if it has the $(k-1)$-skeleton of a cube.  A
cubical $d$-sphere is \emph{(cubically) neighborly} if it is cubically $\lfloor(d+1)/2\rfloor$-neighborly. This is
similar to the simplicial case and not as in Babson, Billera, and Chan, where a $k$-neighborly cubical sphere has the
$k$-skeleton of a cube. In the following $[k]$ denotes the set of positive integers $\{1,\ldots,k\}$.

The next two observations are also implicit in the work of Shemer~\cite{MR693351}.  We start off with a lemma connecting
the neighborliness of the cone over a simplicial ball with the neighborliness of the simplicial ball.

\begin{lem}\label{lem:NeighborlyBBCSequences}
  Let $B$ be a simplicial $(d-1)$-ball on the vertex set $[i-1]$, and let $T = B * i$ its $d$-dimensional cone for some
  $d \geq 2$. Then $\partial T$ is $k$-neighborly if and only if $\partial B$ is $(k-1)$-neighborly and $B$ is
  $k$-neighborly.
\end{lem}
\begin{proof}
  The boundary $\partial T$ is equal to $(\partial B * i) \cup B$ since $T = B * i$. If $\partial T$ is $k$-neighborly
  all $k$-subsets of $[i]$ are faces of $\partial T$.  Because $\partial T$ is the union of $\partial B * i$ and $B$,
  all $(k-1)$-subsets of $[i-1]$ must be contained in $\partial B$ and $B$ must contain all $k$-subsets of $[i-1]$.
  This means that $\partial B$ is $(k-1)$-neighborly and $B$ is $k$-neighborly.
  
  Conversely, if $\partial B$ is $(k-1)$-neighborly and $B$ is $k$-neighborly then $\partial T = (\partial B * v) \cup
  B$ obviously contains all the $k$-subsets of $[i]$.
\end{proof}

We call a BBC sequence $(T_i)_{i=d+1}^{n}$ \emph{neighborly} if the final boundary sphere $\partial T_n$ is neighborly.
With Lemma~\ref{lem:NeighborlyBBCSequences} and induction we obtain a characterization of neighborly BBC sequences.

\begin{prop}
  Let $(T_i)_{i=d+1}^{n}$ be a BBC sequence of $d$-balls for $d \geq 2$ with $B_i$ defined as in
  Definition~\ref{defn:BBCSequences}.  Then the following are equivalent:
  \begin{enumerate}
  \item \label{item:1} $(T_i)_{i=d+1}^{n}$ is a neighborly BBC sequence.
  \item \label{item:2} $\partial T_i$ is a neighborly $(d-1)$-sphere for all $i = d+1,\ldots,n$.
  \item \label{item:3} For all $i=d+1,\ldots,n$ the ball $B_{i-1}$ is $\lfloor d/2 \rfloor$-neighborly, and the
    sphere $\partial B_{i-1}$ is $(\lfloor d/2 \rfloor-1)$-neighborly.
  \end{enumerate}
\end{prop}

Neighborly BBC sequences arise naturally from neighborly simplicial polytopes as in Remark~\ref{rem:BBCSequences}. So
the above definition is a generalization of the sequences of pulling triangulations of cyclic polytopes originally used
by Babson, Billera, and Chan.

\begin{cor}\label{cor:NSPandINS}
  Take an arbitrary ordering of the vertices of a neighborly simplicial polytope. Then there exists a realization such
  that the induced pulling triangulations form a neighborly BBC sequence.
\end{cor}

The neighborliness of a simplicial complex is preserved by mirroring, that is, the mirror complex of $k$-neighborly
simplicial complex is a $(k+1)$-neighborly cubical complex: In its sign vector representation a $k$-neighborly simplicial
complex $\Delta$ on $n$ vertices contains all vectors in $\{ 0,1 \}^n$ with $k$ zeros, i.e. all $k$-subsets of $[n]$.
Thus its mirror complex $\mirror{\Delta}$ contains all vectors in $\{ 0, \pm1 \}^n$ with $k$ zeros. These vectors
represent the $k$-skeleton of the $n$-cube, which means $\mirror{\Delta}$ is $(k+1)$-neighborly.

Cubical neighborliness is also preserved by certain fissuring operations. Given a $k$-neighborly pure cubical complex
$C$ and a subcomplex $S$ whose boundary is $k$-neighborly, then the cubical fissure between $S$ and its complement in
$C$ is again a $k$-neighborly cubical complex.

This yields neighborly cubical spheres from any neighborly BBC sequence with the same construction as in
Theorem~\ref{thm:BBC}.

\begin{cor}
  \label{cor:NeighborlyBBC}
  Let $\cT$ be a neighborly BBC sequence of simplicial $(d-1)$-balls. Then the cubical complexes $\bbc{\cT}$ is a
  neighborly cubical $d$-sphere.
\end{cor}

\section{Cubical Spheres from BBC sequences}

Next we will derive a sign vector representation for the cubical sphere constructed from a BBC sequence.  In
Section~\ref{sec:NeighborlyCubicalSpheres} we show that the neighborly cubical spheres build from special vertex
orderings of cyclic polytopes are indeed isomorphic to the neighborly cubical polytopes studied in~\cite{MR1758054}.
Further we construct a non-polytopal neighborly cubical sphere based on a BBC sequence obtained from the non-polytopal
Altshuler $3$-sphere on $10$ vertices (cf. Altshuler~\cite{MR0467531} and Bokowski and Garms~\cite{MR919873}).

\subsection{A Combinatorial Description}

The following theorem states a purely combinatorial description of the cubical spheres constructed from BBC sequences.
Its characterization as a subcomplex of a high dimensional cube is close to the Cubical Gale Evenness
Condition~\cite[Theorem~18]{MR1758054}.

\begin{thm}
  \label{thm:GeneralizedBBC}
  Let $\cT=(T_i)_{i=d}^{n-1}$ be a BBC sequence of simplicial $(d-1)$-balls with $n>d>2$.  Then the facets of the
  cubical $d$-sphere $\bbc{\cT}$ correspond to the following list of $d$-faces of the $n$-dimensional cube $\cube{n}$;
  they are represented as sign vectors $\alpha \in \{0,\pm1\}^n$ with exactly $d$ zeros.  The type $t$ of a facet
  corresponds to the number of trailing non-zero entries in $\alpha$:
  \begin{itemize}
  \item \textbf{type $t=0$}: $\alpha_n = 0$ and $|\alpha^{(n-1)}| :=
    (|\alpha_1|,\ldots,|\alpha_{n-1}|) \in \partial T_{n-1}$
  \item \textbf{type $0 < t < n-d$}:  $\alpha = (\alpha^{(n-t-1)},0,\alpha_{n-t+1}=\sigma,-1,\ldots,-1)$, where
    \\
    $\sigma \in \{\pm1\}$, $\alpha^{(n-p-1)} \in \{0,\pm 1\}^{n-t-1}$ with
    \begin{enumerate}
    \item $|\alpha^{(n-t-1)}|\in \partial T_{n-t-1}$ and
    \item if $\sigma = +1$, then $|(\alpha^{(n-t-1)},0)| \in T_{n-t}$;\\
      if $\sigma = -1$, then $|(\alpha^{(n-t-1)},0)| \not\in T_{n-t}$.
    \end{enumerate}
  \item \textbf{type $t = n-d$}: $\alpha =
    (0,\ldots,0,\sigma,-1,\ldots,-1)$ with $\sigma \in \{-1, +1\}$.
  \end{itemize}
\end{thm}

\begin{proof}
  We prove that the vector representation of the facets of the cubical sphere given by the theorem corresponds to the
  facets of the inductive definition of $S_k$ used in Theorem~\ref{thm:BBC}: 
  \begin{align*}
    S_{k+1} & = \left(\mirror{T_k} \times \{+1\}\right) \cup \left(\partial\mirror{T_k} \times \{0\}\right) \cup
    \left((S_k \setminus \mirror{T_k}) \times \{-1\}\right)
  \end{align*}  
  It was already shown in Theorem~\ref{thm:BBC} that all $S_k$ are cubical spheres.  We proceed by induction on $k$.
  The sphere $S_{d+1}$ is the boundary of the $(d+1)$-cube for $k = d+1$. The facets of $S_{d+1}$ are the vectors
  $\alpha \in \{ 0, \pm1\}^{d+1}$ with exactly one non-zero entry. The facets of $T_d = \simplex{d-1}$ are the vectors
  in $\{0,1\}^d$ with one $1$ as well.  Thus all the facets of $S_{d+1}$ are either of type $0$ or type $1$ shown in
  Table~\ref{tab:FacetsSd1}.
  \begin{table}[tb]
    \caption{The facets of $S_{d+1}$ are the facets of the $(d+1)$-cube.\label{tab:FacetsSd1}}
    \begin{center}
      \renewcommand{\arraystretch}{1.5}
      \begin{tabular*}{.4\linewidth}{@{\extracolsep{\fill}}c@{\qquad}cccc@{}}\toprule
        Type & \multicolumn{4}{c}{Facets}\\ \midrule
        $0$ & $\pm $     &  $0$       & $\cdots$ & $0$        \\
        $0$ & $0$        & $\pm $     &      $0$ & $\vdots$ \\
        $0$ & $\vdots$   & $\ddots$   & $\ddots$ & $0$        \\
        $1$ & $0$        & $\cdots$   &      $0$ & $\pm $   \\ \bottomrule
      \end{tabular*}
    \end{center}
  \end{table}
  
  The inductive step is split into two claims showing the vector description and the inductive definition of $S_k$ yield
  the same combinatorics.

  \begin{claim}
    Each facet in the vector representation given by the theorem is also a facet in the inductive definition via
    fissuring.
  \end{claim}
  
  We analyze all types of facets of the vector representation.
  \begin{itemize}
  \item \textbf{type $0$:} The facets of type $0$ are the vectors $\alpha = (\alpha^{(k)},0) \in
    \{0,\pm1\}^{k+1}$ with $|\alpha^{(k)}| \in \partial T_k$. This is equivalent to $\alpha^{(k)} \in \mirror{\partial
      T_k}$ and thus $(\alpha^{(k)},0) \in \mirror{\partial T_k} \times \{0\}=\partial\mirror{T_k}\times\{0\}$.
  \item \textbf{type $1$:} The facets of type $1$ are the vectors $\alpha = (\alpha^{(k-1)},0,\sigma)$ with
    $|\alpha^{(k-1)}| \in \partial T_{k-1}$ and
    \begin{itemize}
    \item $\sigma = +1$ if $(\alpha^{(k-1)},0)\in T_k$ or
    \item $\sigma = -1$ if $(\alpha^{(k-1)},0) \not\in T_k$.
    \end{itemize}
    By induction $(\alpha^{(k-1)},0)$ is a type $0$ facet of $S_k$.  Since $T_k = B_{k-1} * v_k$ with $B_{k-1} \subseteq
    \partial T_{k-1}$ we obtain:
    \begin{itemize}
    \item if $(|\alpha^{(k-1)}|,0) \in T_k$ then $(\alpha^{(k-1)},0) \in \mirror{T_k}$ and thus\\ $(\alpha^{(k-1)},0,+1)
      \in S_{k+1}$;
    \item if $(|\alpha^{(k-1)}|,0) \not\in T_k$ then $(\alpha^{(k-1)},0) \in S_k \setminus \mirror{T_k}$ and thus\\
      $(\alpha^{(k-1)},0,-1) \in S_{k+1}$.
    \end{itemize}
  \item \textbf{type $t=2,\ldots,k-d$:} The facets of type $t$ are the vectors
    \[
    (\alpha^{(k-t)},0,\sigma,-1,\ldots,-1)\in \{0,\pm1\}^{k+1}
    \]
    with $\sigma\in\{\pm1\}$ and $\alpha^{(k-t)} \in \{0,\pm1\}^{k-t}$ are the facets of type $t$. Taking only the first
    $k$ entries of $\alpha = (\alpha^{(k)},-1)$, it follows by induction that $\alpha^{(k)}$ is a type $t-1$ facet of
    $S_k$.  Since all facets of $\mirror{T_k}$ contain the vertex $k$, the face $\alpha^{(k)}$ is not contained in
    $\mirror{T_k}$. Further it is not contained in $\mirror{\partial T_k}$ because it has $d$ zero entries. Thus
    $\alpha^{(k)}$ is a facet of $S_{k} \setminus \mirror{T_k}$ and $(\alpha^{(k)},-1) \in S_{k+1}$.
  \end{itemize}
  
  \begin{claim}
    Each facet in the inductive definition via fissuring is also a facet in the vector representation given by the
    theorem.
  \end{claim}
  
  There are three different kinds of facets of $S_{k+1}$ according to the inductive definition. We will determine the
  type of each kind of facet.
  \begin{itemize}
  \item The facets of $\mirror{T_{k}}\times\{+1\}$ are the vectors $(\alpha^{(k)},+1) \in \{0,\pm1\}^{k+1}$ with
    $|\alpha^{(k)}| \in T_k$ and $\alpha^{(k)}_k = 0$ since by definition of the BBC sequence all the facets of $T_k$
    contain the vertex $k$.  These are some of the facets of type $1$ of $S_{k+1}$.
  \item The facets of $\mirror{\partial T_k}\times\{0\}$ are the vectors $(\alpha^{(k)},0) \in \{0,\pm1\}^{k+1}$ with
    $|\alpha^{(k)}| \in \partial T_k$. These are the facets of type $0$ of $S_{k+1}$.
  \item Let $\alpha^{(k)} \in S_{k}\setminus \mirror{T_k}$. By induction, we get the vector representation of the facets
    of $S_k$.
    
    If $\alpha^{(k)}$ is a facet of type $t$ with $0<t\leq k-d$ then $\alpha_{k}^{(k)} \neq 0$ and thus is not in
    $\mirror{T_k}$.  By appending a $-1$ we get a facet of type $t+1$ of $S_{k+1}$.
    
    The facets $(\alpha^{(k-1)},0) \in S_k$ with $|\alpha^{(k-1)}| \in \partial T_{k-1}$ are of type $0$. Those facets
    may be partitioned into two parts: $(|\alpha^{(k-1)}|,0) \in T_k$ or $(|\alpha^{(k-1)}|,0)\not\in T_k$.  If
    $(|\alpha^{(k-1)}|,0) \in T_k$ then $(\alpha^{(k-1)},0) \in \mirror{T_k}$ and thus $(\alpha^{(k-1)},0,+1)\in
    S_{k+1}$.  Similarly, if $(|\alpha^{(k-1)}|,0) \not\in T_k$ then $(\alpha^{(k-1)},0) \in S_k\setminus\mirror{T_k}$
    and thus $(\alpha^{(k-1)},0,-1)\in S_{k+1}$.  These are the facets of type~$1$ of~$S_{k+1}$. 
  \end{itemize}
  
  Hence the cubical sphere $\bbc{(T_i)_{i=d}^{n-1}} = S_{n}$ has the facets given by the theorem.
\end{proof}

If $\cT=(T_i)_{i=d}^{n-1}$ is a neighborly BBC sequence, $n>d>2$, then each simplicial $(d-2)$-sphere $\partial T_i$ is
neighborly.  If additionally $d$ is odd then $s(i,d-2)$, the number of facets of $\partial T_i$, is determined by the Dehn-Sommerville
equations, see Gr\"unbaum~\cite[\S9.2]{MR1976856}.  We have
\[ s(i,d-2)=\binom{i-\frac{d-1}{2}}{i+1-d}+\binom{i-\frac{d+1}{2}}{i+1-d}=\frac{2i}{2i-d+1}\binom{i-\frac{d-1}{2}}{i+1-d}. \]
Adding up this results in a formula for $f(n,d)$, the number of facets of the neighborly cubical sphere $\bbc{\cT}$.
The same formula occurs in~\cite[Corollary~19]{MR1758054} where, however, then even and the odd-dimensional case are
erroneously exchanged:

\begin{cor}\label{cor:n_facets}
  Let $\cT=(T_i)_{i=d}^{n-1}$ be a neighborly BBC sequence of simplicial $(d-1)$-balls with $n>d>2$ and $d$ odd.  Then
  the number of facets of $\bbc{\cT}$ is given by
  \begin{eqnarray*}
    f(n,d) &=& 2(d+1) + \sum_{t=0}^{n-d-2} s(n-t-1,d-2)2^{n-t-d}\\
           &=& 2(d+1) + \sum_{k=d+2}^{n} \frac{2k-2}{2k-1-d}\binom{k-\frac{d+1}{2}}{k-d}2^{k-d}.
  \end{eqnarray*}
\end{cor}

\begin{exmp}\label{exmp:BBC}
  A neighborly cubical $3$-sphere with the graph of the $6$-cube may be constructed from a BBC sequence of a pentagon.
  Its $f$-vector is $(64,192,192,64)$.  The neighborly BBC sequence contains three elements, the triangle $T_3$, and the
  pulling triangulations of the $4$- and the $5$-gon, $T_4$ and $T_5$, respectively. If the vertices are given in cyclic
  order then the vector representations of $T_i$ and $\partial T_i$ for $i=3,4,5$ are given in
  Table~\ref{tab:exmp:BBC:seq}.
  \begin{table}[tb]
    \caption{BBC sequence of a pentagon as described in Example~\ref{exmp:BBC}. The vectors in $\{0,1\}^i$ are the
      characteristic functions of the complement of the vertex sets in $[i]$, for example $11001$ represents
      $\{3,4\}$. The subcomplex $B_{i-1}$ of $\partial T_{i-1}$ such that $T_i$ is the join of $B_{i-1}$ with the vertex
      $i$ are printed in bold.
      \label{tab:exmp:BBC:seq}}
    \begin{center}
      \renewcommand{\arraystretch}{0.9}
      \begin{tabular*}{\linewidth}{@{\extracolsep{\fill}}l@{\hspace{8.5ex}}lll@{\hspace{6.5ex}}lllll@{}}\toprule
              & \multicolumn{3}{c}{Facets of $T_i$} & \multicolumn{5}{c}{Facets of $\partial T_i$}\\ \midrule
        $T_3$ & $000$   &         &         & $\mathbf{001}$   & $\mathbf{100}$   & $010$   &         &\\ \addlinespace
        $T_4$ & $0010$  & $1000$  &         & $\mathbf{0011}$  & $\mathbf{1001}$  & $\mathbf{1100}$  & $0110$  &\\ \addlinespace
        $T_5$ & $00110$ & $10010$ & $11000$ & $00111$ & $10011$ & $11001$ & $11001$ & $01110$\\ \bottomrule      
      \end{tabular*}
    \end{center}
  \end{table}
  The mirror complexes of these are obtained by simply replacing the $1$'s by $\pm1$'s. So according to the inductive
  definition we start with $S_4$, the boundary of the $4$-cube and fissure twice along $\partial T_4$ and then along
  $\partial T_5$. This yields the facet description of $S_6$ listed in Table~\ref{tab:exmp:BBC:facets}.
  \begin{table}[bt]
    \caption{Facets of the neighborly cubical $3$-sphere $S_6$ sorted by type; see Example~\ref{exmp:BBC}. The facets
      occur in different multiplicities depending on the number of $\pm$ entries in the corresponding row: For example,
      every row of type $0$ represents eight facets.
      The rows of the $4\times4$ minor printed in bold are the facets of $S_4$, that is, the boundary of the $4$-cube. 
      The facets of $S_5$ are the $8\times5$ minor printed in bold and plain.
      The horizontal lines underline the inductive structure.
      \label{tab:exmp:BBC:facets}}
    \begin{center}
      \renewcommand{\arraystretch}{0.9}
      \begin{tabular*}{.75\linewidth}{@{\extracolsep{\fill}}c@{\qquad}rrrrrr@{}}\toprule
        Type & \multicolumn{6}{c}{Facets}\\ \midrule
        $3$  & $\mathbf{0}$ & $\mathbf{0}$ & $\mathbf{0}$ & $\mathbf{\pm }$ & $- $ & $\mathbf{- }$\\ \addlinespace
        $2$  & $\mathbf{0}$ & $\mathbf{0}$ & $\mathbf{\pm }$ & $\mathbf{0}$ & $+ $ & $\mathbf{- }$\\
        $2$  & $\mathbf{0}$ & $\mathbf{\pm }$ & $\mathbf{0}$ & $\mathbf{0}$ & $- $ & $\mathbf{- }$\\
        $2$  & $\mathbf{\pm }$ & $\mathbf{0}$ & $\mathbf{0}$ & $\mathbf{0}$ & $+ $ & $\mathbf{- }$\\ \cmidrule{2-5} \addlinespace
        $1$  & $0$ & $0$ & $\pm $ & $\pm $ & $0$ & $\mathbf{+ }$\\
        $1$  & $\pm $ & $0$ & $0$ & $\pm $ & $0$ & $\mathbf{+ }$\\
        $1$  & $\pm $ & $\pm $ & $0$ & $0$ & $0$ & $\mathbf{+ }$\\
        $1$  & $0$ & $\pm $ & $\pm $ & $0$ & $0$ & $\mathbf{- }$\\ \cmidrule{2-6}  \addlinespace
        $0$ & $\mathbf{0}$ & $\mathbf{0}$ & $\mathbf{\pm }$ & $\mathbf{\pm }$ & $\mathbf{\pm }$ & $\mathbf{0}$\\
        $0$ & $\mathbf{\pm }$ & $\mathbf{0}$ & $\mathbf{0}$ & $\mathbf{\pm }$ & $\mathbf{\pm }$ & $\mathbf{0}$\\
        $0$ & $\mathbf{\pm }$ & $\mathbf{\pm }$ & $\mathbf{0}$ & $\mathbf{0}$ & $\mathbf{\pm }$ & $\mathbf{0}$\\
        $0$ & $\mathbf{\pm }$ & $\mathbf{\pm }$ & $\mathbf{\pm }$ & $\mathbf{0}$ & $\mathbf{0}$ & $\mathbf{0}$\\
        $0$ & $\mathbf{0}$ & $\mathbf{\pm }$ & $\mathbf{\pm }$ & $\mathbf{\pm }$ & $\mathbf{0}$ & $\mathbf{0}$\\ \bottomrule
      \end{tabular*}
    \end{center}
  \end{table}
\end{exmp}

\subsection{Neighborly Cubical Polytopes}
\label{sec:NeighborlyCubicalSpheres}
In this section we show that for very particular neighborly BBC sequences, the neighborly cubical spheres constructed
therefrom are isomorphic to the boundaries of the neighborly cubical polytopes described in~\cite[Theorem~18]{MR1758054}
and thus polytopal. The neighborly cubical polytopes have the following sign vector representation.

\begin{thm}[Cubical Gale Evenness Condition]
\label{thm:CubicalGale}
The facets of the neighborly cubical polytope $\ncp{d+1}{n}$ are given by vectors $\alpha \in \{0, \pm 1\}^n$ with $d$
zeros. They are classified by the number $t$ of leading $\pm 1$'s:
  \begin{itemize}
  \item \textbf{type $t=0$}: $\alpha_1 = 0$, and $|\alpha|$ satisfies the simplicial Gale Evenness Condition: between
    any two values $\alpha_i, \alpha_j \in \{\pm 1\}$ there is an even number of zeros.
  \item \textbf{type $0 < t < n-d$}:  $\alpha = (-1,+1,\ldots,(-1)^{t-1},\sigma,0,\alpha^{(n-t-1)})$,
    with\\
    $\sigma \in \{\pm1\}$, $\alpha^{(n-t-1)} \in \{0,\pm 1\}^{n-t-1}$ and:
    \begin{enumerate}
    \item $|\alpha^{(n-t-1)}|$ satisfies the simplicial Gale Evenness
      Condition, and
    \item if $\sigma = (-1)^{t+1}$, then $\alpha^{(n-t-1)}$ starts with an
      even number of zeros; \\
      if $\sigma = (-1)^{t}$, then $\alpha^{(n-t-1)}$ starts with an odd number of zeros.
    \end{enumerate}
  \item \textbf{type $t = n-d$}: $\alpha = (-1,+1,\ldots,(-1)^{t-1},\sigma,0,\ldots,0)$ with $\sigma \in \{-1, +1\}$.
  \end{itemize}
\end{thm}

Consider the cyclic $(d-1)$-polytope on $n$ vertices $\cyclic{d-1}{n}$ given as the convex hull of the vertices $v_j =
(j,j^2, \ldots, j^d)$ for $j = 1, \ldots, n$. The facets of the cyclic polytope are given by Gale's Evenness
Condition, cf.~\cite{MR0152944} and~\cite[Theorem~0.7]{MR1311028}. Let $T_i$ denote the pulling triangulation of
$\cyclic{d-1}{i}$ with respect to the vertex $v_i$. The facets of $T_i$ according to our vector notation for simplicial
complexes are the vectors $\phi \in \{0,1\}^i$ with $d$ zeros, such that:
\begin{enumerate}
\item $\phi_i = 0$,
\item $(\phi_1,\ldots,\phi_{i-1})$ satisfies Gale's Evenness Condition, and
\item $\phi$ ends with an odd number of zeros.
\end{enumerate}
The facets of the neighborly cubical sphere constructed from the corresponding neighborly BBC sequence
$(T_i)_{i=d}^{n-1}$ are readily derived from Theorem~\ref{thm:GeneralizedBBC}. These specific spheres are polytopal,
since they are isomorphic to the neighborly cubical polytopes.

\begin{cor}
  \label{cor:NCSeqNCP}
  Let $d \geq 3$ and $T_i$ be the pulling triangulation of the cyclic polytope $\cyclic{d-1}{i}$ with respect to the
  last vertex as above. Then the neighborly cubical sphere $\bbc{(T_i)_{i=d}^{n-1}}$ and the boundary of the neighborly
  cubical polytope $\ncp{d+1}{n}$ are combinatorially isomorphic. The isomorphism is given by inverting the order and
  then flipping the even bits:
  \begin{align*}
  \Phi : \{0, \pm 1\}^n &\to \{0, \pm 1\}^n, \\
  (\alpha_1,\alpha_2,\ldots,\alpha_n) &\mapsto 
  (\alpha_n,-\alpha_{n-1},\ldots,(-1)^{n}\alpha_2,(-1)^{n+1}\alpha_1)
  \end{align*}
\end{cor}

\begin{rem}
  We briefly summarize some results studied in~\cite[Section~3.6.3]{DiplomNCS}. The neighborly cubical $d$-spheres
  constructed from BBC sequences of cyclic $(d-1)$-polytopes on $d+1$ vertices are combinatorially isomorphic.
  
  Since in \textbf{even} dimension the automorphism group of every cyclic polytope acts transitively on its vertices,
  its pulling triangulations are combinatorially isomorphic. Thus the neighborly cubical spheres constructed from
  arbitrary vertex orderings of even dimensional cyclic polytopes are combinatorially unique and thus all polytopal.
  This includes the Example~\ref{exmp:BBC}.
  
  In \textbf{odd} dimension $d-1$, the automorphism group of a cyclic polytope with more than $d+1$ vertices does not
  act transitively on the vertices. Thus we were able to construct non-isomorphic neighborly cubical $4$-spheres from
  different triangulation of the cyclic $3$-polytope on six vertices.
\end{rem}

\subsection{A Non-polytopal Neighborly Cubical Sphere}
\label{sec:NonPolytopalSphere}

To algorithmically decide the polytopality of spheres is generally possible but of outrageous complexity, see Bokowski
and Sturmfels~\cite{MR1009366}, Richter-Gebert~\cite{MR1482230}, and Basu et al.~\cite{MR1998147}.  If standard
heuristic methods fail to decide the polytopality of a given sphere typically a suitable ad-hoc proof is the only
remedy.  A good example is the $3$-sphere $\mathcal{M}^{10}_{425}$ found by Altshuler in an exhaustive enumeration of
all neighborly simplicial $3$-manifolds on ten vertices~\cite{MR0467531}.  Its polytopality at first could not be
decided, and so it was up to Bokowski and Garms~\cite{MR919873} to prove that $\mathcal{M}^{10}_{425}$ does not admit
any convex realization.  Below we take this very example as the starting point for a construction of a neighborly
cubical sphere which cannot be polytopal in view of the following result.

\begin{thm}\label{thm:polytopal}
 Let $\cT=(T_i)_{i=d}^{n-1}$ be a BBC sequence of simplicial $(d-1)$-balls, neighborly or not, such that the cubical
 $d$-sphere $\bbc{\cT}$ is polytopal.  Then the simplicial $(d-2)$-sphere $\partial T_{n-1}$ necessarily is polytopal.
\end{thm}

\begin{proof}
  Suppose that $\bbc{\cT}$ is isomorphic to the boundary complex of a convex cubical $(d+1)$-polytope
  $P\subset\RR^{d+1}$.  Consider the edge $e$ of~$P$ corresponding to the sign vector $(+1,\dots,+1,0)$ of length~$n$.
  Choose an affine hyperplane~$H$, parallel to~$e$, which separates the edge~$e$ from the $2^n-2$ vertices of~$P$ not
  contained in~$e$.  Next we choose a second affine hyperplane $H'$, orthogonal to~$e$, which separates the two vertices
  of~$e$.  Then, since $H$ and $H'$ are not parallel to each other, $P/e=P\cap H\cap H'$ is a $(d-1)$-dimensional convex
  polytope, the \emph{edge figure} of $e$ with respect to~$P$.  Its face lattice is isomorphic to the filter of $e$,
  that is, the faces of~$P$ containing~$e$, in the face lattice of~$P$.
  
  Now the facets of~$P$ which contain~$e$ are exactly the facets of type~$0$ without negative signs in their sign vector
  representation.  From Theorem~\ref{thm:GeneralizedBBC} we conclude that $\partial (P/e)$ is isomorphic to $\partial
  T_{n-1}$, and hence the claim.
\end{proof}

As a consequence, each BBC sequence of simplicial $(d-1)$-balls with the property that its final boundary sphere is
non-polytopal yields a non-polytopal cubical sphere.  Note, however, that in contrast to polytopal spheres, see
Corollary~\ref{cor:NSPandINS}, for non-polytopal spheres there is no standard procedure to obtain a corresponding BBC
sequence.  Moreover, it follows from work of Altshuler that there is a simplicial $3$-sphere on ten vertices which does
not admit any BBC sequence~\cite{MR0397744}.

\begin{thm}\label{thm:Altshuler}
  There is a non-polytopal neighborly cubical $5$-sphere with $2^{11}=2048$ vertices and $f(11,5)=3584$ facets.  Its
  complete $f$-vector is
  \[ f=(2048,11264,28160,33280,17920,3584). \]
\end{thm}
\begin{proof}
  Table~\ref{tab:BBCSequenceOfAlthulerSphere} lists a neighborly BBC sequence $\cA=(A_i)_{i=5}^{10}$ of simplicial
  $4$-balls with the property that the boundary of the final ball $A_{10}$ is isomorphic to the Altshuler $3$-sphere
  $\mathcal{M}^{10}_{425}$.
  
  The simplicial $4$-ball $A_8$ is a pulling triangulation of the neighborly simplicial $4$-polytope on eight vertices
  which occurs as $P_{36}^8$ in the list of Gr\"unbaum and Sreedharan~\cite{MR0215182}; this is one of the two
  non-cyclic neighborly simplicial polytopes with these parameters.
  
  The simplicial $4$-ball $A_9$ is a pulling triangulation of a neighborly simplicial $4$-polytope on nine vertices: The
  previous simplicial $3$-sphere $\partial A_8$ is separated by the $2$-sphere $\partial B_8$ into the $3$-balls $B_8$
  and its complement $B_8'$.  The boundary $\partial A_9$ now is directly obtained from $\partial A_8$ by first removing
  $B_8'$ and then inserting the cone over $\partial B_8$ with the new vertex~$9$ as the apex.  Since $\partial A_8$ is
  polytopal and $B_8'$ is a $3$-simplex which is stacked once over each of its four facets it follows that $A_9$ is
  again polytopal: The vertex~$9$ can be chosen in a way such that it is \emph{exactly beyond} all the facets of~$B_8'$
  in the sense of Shemer~\cite[page~301]{MR693351}.
  
  The sphere $\partial A_{10}$ is directly obtained from $\partial A_9$ and isomorphic to~$\mathcal{M}^{10}_{425}$ and
  thus not polytopal due to Bokowski and Garms~\cite{MR919873}.  From Theorem~\ref{thm:polytopal} it follows that the
  neighborly cubical sphere $\bbc{\cT}$ is not polytopal.
\end{proof}

\begin{figure}[htbp]
  \centering
  \subfigure[Bipyramid over triangle: $\partial B_5$.]{
    \begin{overpic}[width=.3\textwidth]{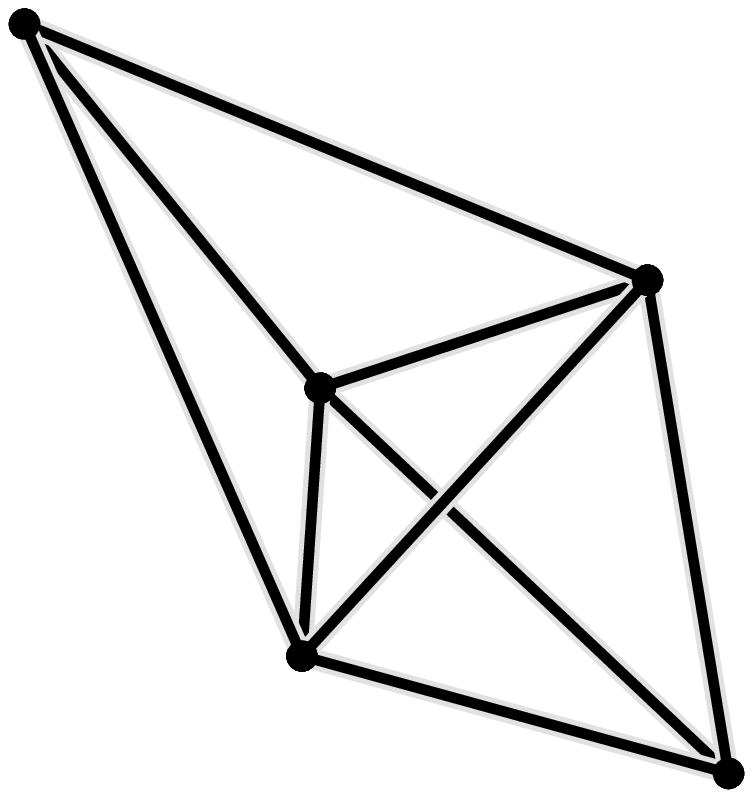}
      \put(80,68){$0$}
      \put(5,93){$1$}
      \put(90,10){$2$}
      \put(44,60){$3$}
      \put(33,20){$4$}
    \end{overpic}
  }
  \subfigure[$3$-simplex stacked on two facets: $\partial B_6$.]{
    \begin{overpic}[width=.3\textwidth]{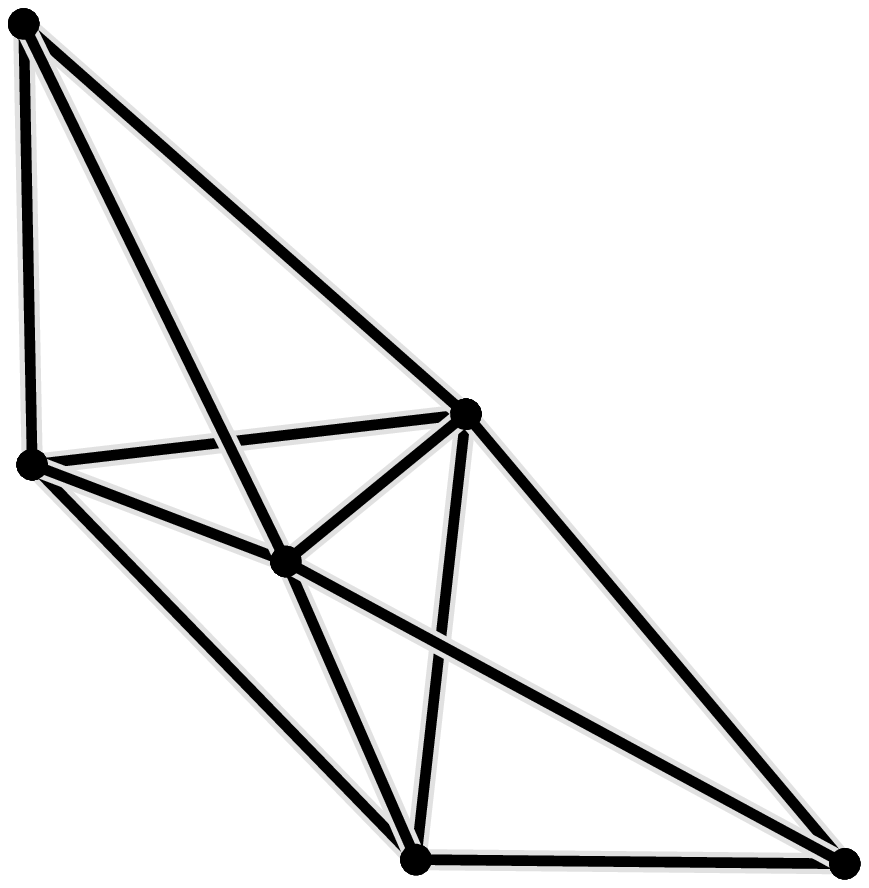}
      \put(35,35){$0$}
      \put(102,10){$1$}
      \put(63,60){$2$}
      \put(6,97){$3$}
      \put(45,5){$4$}
      \put(7,50){$5$}
    \end{overpic}
  }
  \subfigure[Cone over a strip of four triangles: $\partial B_7$.]{
    \begin{overpic}[width=.3\textwidth]{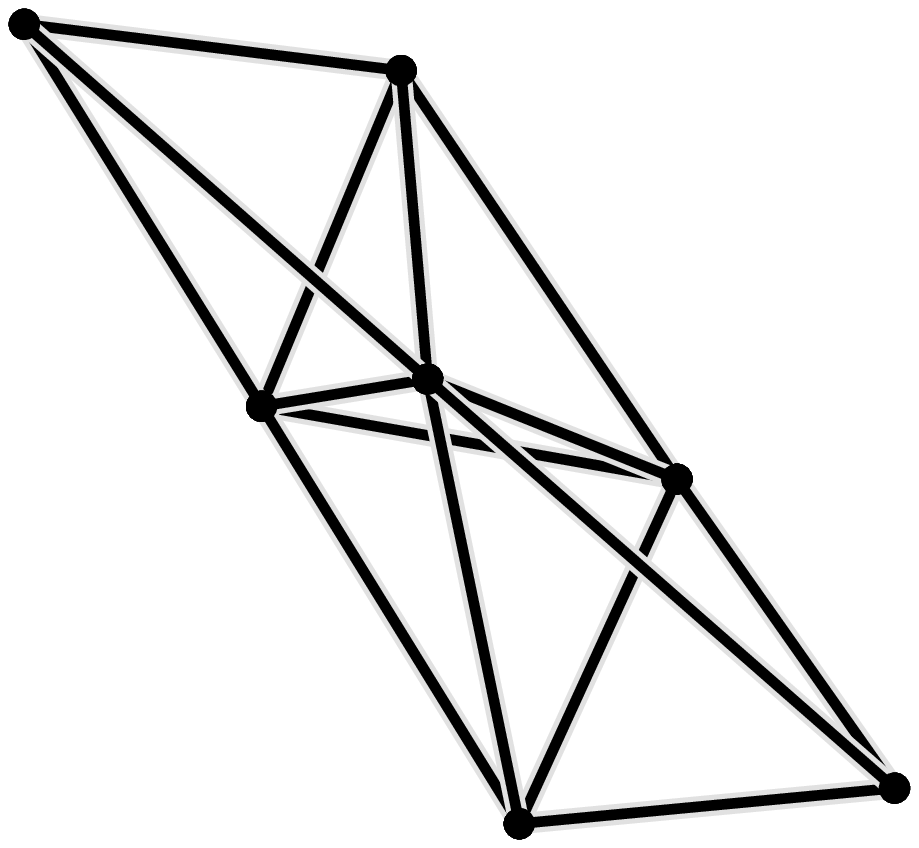}
      \put(-2,90){$0$}
      \put(50,87){$1$}
      \put(78,45){$2$}
      \put(100,10){$3$}
      \put(22,48){$4$}
      \put(48,4){$5$}
      \put(51,58){$6$}
    \end{overpic}
  }
  \subfigure[$3$-simplex stacked on all its facets: $\partial B_8$.]{
    \begin{overpic}[width=.3\textwidth]{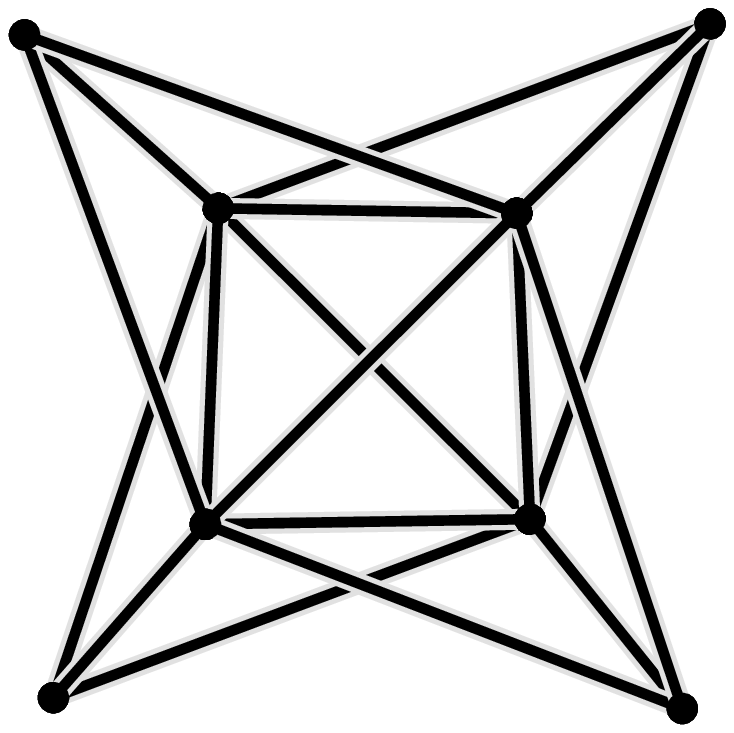}
      \put(10,82){$0$}
      \put(30,60){$1$}
      \put(67,23){$2$}
      \put(92,82){$3$}
      \put(72,60){$4$}
      \put(90,9){$5$}
      \put(12,10){$6$}
      \put(35,23.5){$7$}
    \end{overpic}
  }
  \subfigure[Cone over a strip of six triangles: $\partial B_9$.]{
    \begin{overpic}[width=.3\textwidth]{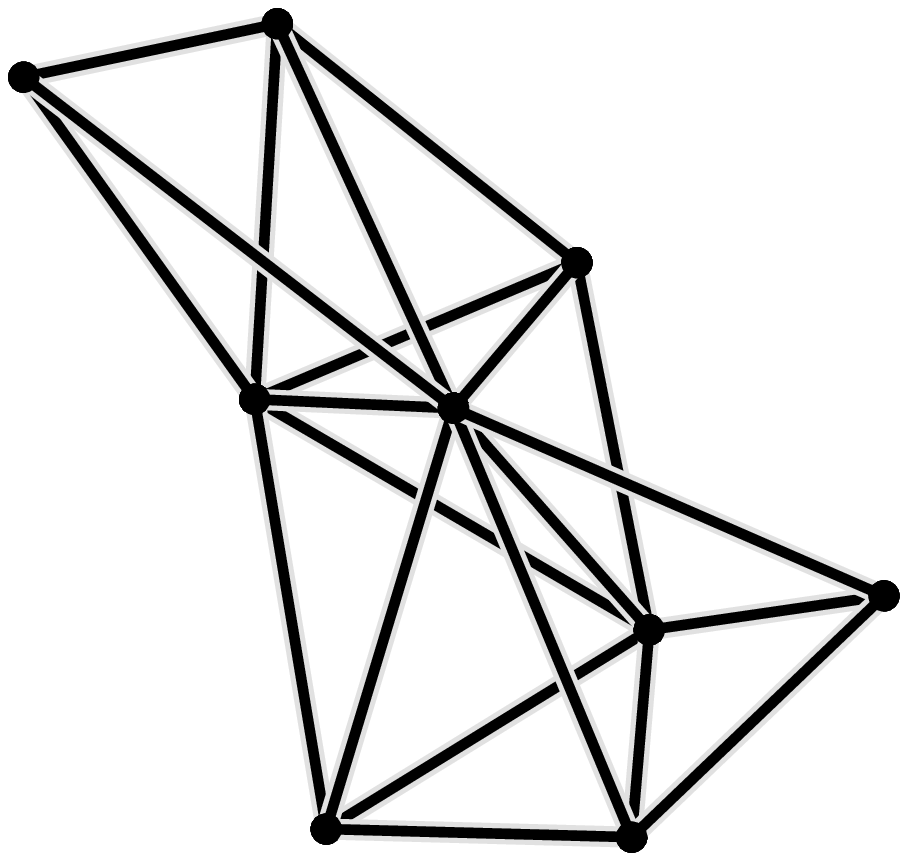}
      \put(68,68){$0$}
      \put(73,32){$1$}
      \put(69,-2){$2$}
      \put(100,30){$3$}
      \put(32,95){$4$}
      \put(0,83){$5$}
      \put(36,-1){$6$}
      \put(23,45){$7$}
      \put(59,52){$8$}
    \end{overpic}
  }
  \caption{The subcomplexes $B_i'=\partial A_i\setminus B_i$ of the BBC sequence of the Altshuler sphere; see
    Theorem~\ref{thm:Altshuler}.  The sphere $\partial B_i=\partial B_i'$ is the link of the new vertex $i$ in the
    sphere $\partial A_{i+1}$.\label{fig:BoundariesOfBks}}
\end{figure}

Figure~\ref{fig:BoundariesOfBks} visualizes the complements of the simplicial $3$-balls $B_5$, $B_6$, $B_7$, $B_8$,
and~$B_9$ and their boundaries.

\begin{rem}
  The polytopal neighborly simplicial $3$-sphere $\partial A_9$ is the unique $3$-sphere on nine vertices with Altshuler
  determinant equal to $103361328$.  It occurs as $N_{19}^9$ in the classification of Altshuler and
  Steinberg~\cite{MR0333981} and as \texttt{manifold\_3\_9\_523} in Lutz~\cite{Lutz_PAGE}.
\end{rem}

\begin{table}[tb]
  \caption{Facets of the neighborly BBC sequence $\cA$ of the non-polytopal Altshuler sphere $\mathcal{M}^{10}_{425}$. The
    simplicial ball $A_{i}$ is the join the vertex $i-1$ and the subcomplex $B_{i-1}$ of the facets of $\partial A_{i-1}$
    printed in bold.}
  \label{tab:BBCSequenceOfAlthulerSphere}
  \begin{center}
    \renewcommand{\arraystretch}{0.9}
    \begin{tabular*}{\linewidth}{@{\extracolsep{\fill}}l@{\qquad}cccc@{\qquad}ccccc}\toprule
               & \multicolumn{4}{c}{Facets of $A_i$} & \multicolumn{4}{c}{Facets of $\partial A_i$} \\ \midrule
      $A_5$    & $ 01234 $&         &         &         &$ \mathbf{0123} $&$ \mathbf{0124} $&$ 0134 $&$ 0234 $&$ \mathbf{1234} $ \\ \addlinespace
      $A_6$    & $ 01235 $&$ 01245 $&$12345  $&$       $&$  \mathbf{0123} $&$ 0124 $&$ \mathbf{0135} $&$ \mathbf{0145} $&$ 0235$ \\
               & $       $&$       $&$       $&$       $&$  0245 $&$ \mathbf{1234} $&$ \mathbf{1345} $&$ \mathbf{2345} $& \\ \addlinespace
      $A_7$    & $01236  $&$ 01356 $&$ 01456 $&$ 12346 $&$  \mathbf{0123} $&$ \mathbf{0126} $&$ \mathbf{0135} $&$ \mathbf{0145} $&$ 0146$ \\
               & $13456  $&$ 23456 $&$       $&$       $&$  \mathbf{0236} $&$ \mathbf{0356} $&$ \mathbf{0456} $&$ \mathbf{1234} $&$ 1246$ \\
               & $       $&$       $&$       $&$       $&$  \mathbf{1345} $&$ \mathbf{2345} $&$ 2356 $&$ 2456 $&$$ \\ \addlinespace

      $A_8$    & $01237 $&$ 01267 $&$ 01357 $&$ 01457 $&$  \mathbf{0123} $&$ \mathbf{0126} $&$ \mathbf{0135} $&$ \mathbf{0145} $&$ 0147 $\\
               & $02367 $&$ 03567 $&$ 04567 $&$ 12347 $&$  \mathbf{0167} $&$ \mathbf{0236} $&$ \mathbf{0356} $&$ \mathbf{0456} $&$ \mathbf{0467} $\\
               & $13457 $&$ 23457 $&$       $&$       $&$  1234 $&$ 1247 $&$ 1267 $&$ \mathbf{1345} $&$ \mathbf{2345} $\\
               & $      $&$       $&$       $&$       $&$  \mathbf{2357} $&$ \mathbf{2367} $&$ 2457 $&$ \mathbf{3567} $&$ \mathbf{4567} $\\ \addlinespace

      $A_9$    & $01238 $&$ 01268 $&$ 01358 $&$ 01458 $&$  \mathbf{0123} $&$ \mathbf{0126} $&$ \mathbf{0135} $&$ \mathbf{0145} $&$ \mathbf{0148} $\\
               & $01678 $&$ 02368 $&$ 03568 $&$ 04568 $&$  \mathbf{0167} $&$ 0178 $&$ \mathbf{0236} $&$ \mathbf{0356} $&$ \mathbf{0456} $\\
               & $04678 $&$ 13458 $&$ 23458 $&$ 23578 $&$  \mathbf{0467} $&$ 0478 $&$ 1238 $&$ 1268 $&$ \mathbf{1345} $\\
               & $23678 $&$ 35678 $&$ 45678 $&$       $&$  \mathbf{1348} $&$ 1678 $&$ \mathbf{2345} $&$ \mathbf{2348} $&$ \mathbf{2357} $\\
               & $      $&$       $&$       $&$       $&$  \mathbf{2367} $&$ \mathbf{2458} $&$ \mathbf{2578} $&$ \mathbf{2678} $&$ \mathbf{3567} $\\
               & $      $&$       $&$       $&$       $&$  \mathbf{4567} $&$ 4578 $&$ $&$ $&$ $\\ \addlinespace

      $A_{10}$ & $01239 $&$ 01269 $&$ 01359 $&$ 01459 $&$  0123 $&$ 0126 $&$ 0135 $&$ 0145 $&$ 0148 $\\
               & $01489 $&$ 01679 $&$ 02369 $&$ 03569 $&$  0167 $&$ 0179 $&$ 0189 $&$ 0236 $&$ 0356 $\\
               & $04569 $&$ 04679 $&$ 13459 $&$ 13489 $&$  0456 $&$ 0467 $&$ 0479 $&$ 0489 $&$ 1239 $\\
               & $23459 $&$ 23489 $&$ 23579 $&$ 23679 $&$  1269 $&$ 1345 $&$ 1348 $&$ 1389 $&$ 1679 $\\
               & $24589 $&$ 25789 $&$ 26789 $&$ 35679 $&$  2345 $&$ 2348 $&$ 2357 $&$ 2367 $&$ 2389 $\\
               & $45679 $&$       $&$       $&$       $&$  2458 $&$ 2578 $&$ 2678 $&$ 2689 $&$ 3567 $\\
               & $      $&$       $&$       $&$       $&$  4567 $&$ 4579 $&$ 4589 $&$ 5789 $&$ 6789 $\\ \bottomrule
    \end{tabular*}
  \end{center}
\end{table}

\section{Polyhedral Surfaces}
\label{sec:PolyhedralSurfaces}

In this section we describe a nice way to realize polyhedral surfaces of `unusually large genus'~\cite{MR727027} in the
Schlegel diagram of neighborly cubical polytopes. The surfaces considered are cubical polyhedral surfaces where each
vertex has degree $q$. They were first described by Coxeter~\cite{C:SkewPolyhedra} in 1937 in terms of reflection
groups.  Ringel~\cite{MR0075586} explicitly described these surfaces whilst analyzing problems concerning the graph of
the $n$-dimensional cube.  He pointed out that the surfaces are of lowest genus among all surfaces on which the graph
of the $n$-cube may be drawn without self intersection. Further he gave an explicit combinatorial description of the
surface as a $2$-dimensional subcomplex of the $n$-cube. McMullen, Schulz, and Wills~\cite{MR657865,MR727027} analyze
equivelar surfaces, a much more general class of polyhedral surfaces which include the cubical surfaces of Coxeter and
Ringel.

\begin{defn}
  An \emph{equivelar surface} $\mathcal{M}_{p,q}$ is a polyhedral surface such that all $2$-faces are $p$-gons and all
  vertices have degree $q$.
\end{defn}

McMullen, Schulz, and Wills were the first to point at the `unusually large genus' of these surfaces.  In particular,
they inductively constructed an embedding of the surfaces $\mathcal{M}_{4,q}$ in $\RR^3$.

\subsection{Equivelar surfaces and mirror complexes}

Let $Q$ be the boundary of a $q$-gon for $q > 2$ and $\mirror{Q}$ its mirror complex. With the vertices of $Q$ labelled
in cyclic order we obtain the same vector representation as Ringel~\cite[page~17]{MR0075586} of $\mirror{Q}$:
\[
\begin{matrix}
   0 & 0 &\pm&\pm&\cdots&\pm&\pm&\pm\\
  \pm& 0 & 0 &\pm&\cdots&\pm&\pm&\pm\\
   . & . & . & . &\cdots& . & . & . \\
   . & . & . & . &\cdots& . & . & . \\
  \pm&\pm&\pm&\pm&\cdots&\pm& 0 & 0 \\
   0 &\pm&\pm&\pm&\cdots&\pm&\pm& 0 
\end{matrix}
\]

The mirror complex of $Q$ is an equivelar surface $\mathcal{M}_{4,q}$ since the vertex figure of every vertex of
$\mirror{Q}$ is isomorphic to $Q$. It is embedded in the $2$-skeleton of the $q$-cube.

The obvious way to realize a $2$-dimensional subcomplex of the $q$-cube in $\RR^5$ is in the Schlegel diagram of the
$6$-dimensional neighborly cubical polytope $\ncp{6}{q}$ of~\cite{MR1758054}: Since the neighborly cubical polytope
$\ncp{6}{q}$ has the $2$-skeleton of the $q$-cube, $\mirror{Q}$ is contained in the boundary of $\ncp{6}{q}$.  The
Schlegel diagram of $\ncp{6}{q}$ is embedded in $\RR^5$, thus $\mirror{Q}$ may be realized in $\RR^5$.

Let $T_i$ be the pulling triangulation of the $i$-gon for $i=3,\ldots, q-1$. Then the facets of type $0$ of
$\mathcal{S}_3(q) = \bbc{(T_i)_{i=3}^{q-1}}$ are a cone over the boundary of the $(q-1)$-gon. Thus the mirror complex of
the $q$-gon and its pulling triangulation are both subcomplexes of $\mathcal{S}_3(q)$. This yields a cubificated
embedding of the surface $\mirror{Q}$ into the neighborly cubical $3$-sphere $\mathcal{S}_3(q)$. By
Corollary~\ref{cor:NCSeqNCP} this sphere is isomorphic to the boundary of the neighborly cubical polytope $\ncp{4}{q}$.
Hence the mirror complex of the $q$-gon can be realized as a subcomplex of the Schlegel diagram of $\ncp{4}{q}$ in
$\RR^3$.  This answers a question of G\"unter M. Ziegler~\cite{ZieglerPC} whether some of the surfaces
from~\cite{MR727027} can be found as subcomplexes of the neighborly cubical polytopes.

The genus of this surface may easily be calculated from its $f$-vector $f(\mirror{Q}) =(2^q,2^{q-1}q,2^{q-2}q)$:
\[ 
g(q) = 1 + 2^{q-3}(q-4) \in \mathcal{O}(f_0\cdot \log f_0).
\]
Thus $q=12$ is the first parameter for which the genus $g(12)=4097$ exceeds the number of vertices, which equals $2^{12}
= 4096$.

\begin{figure}[ht]
  \centering \includegraphics[height=8cm]{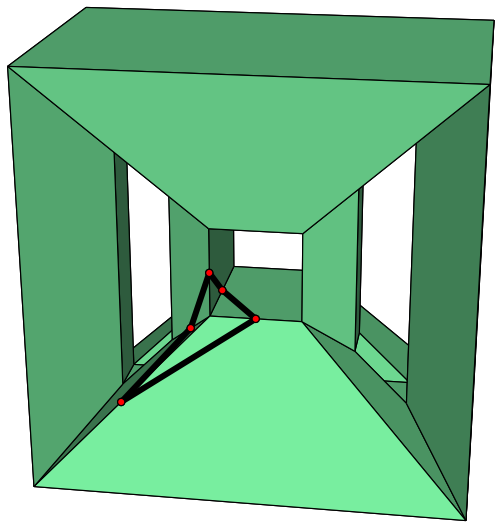}
  \caption{The mirror complex of the $5$-gon is an equivelar surface of type $\mathcal{M}_{4,5}$ of genus $5$. It is
    realized in the Schlegel diagram of the neighborly cubical polytope $\ncp{4}{5}$, which is embedded in $\RR^3$.}
  \label{fig:M_5gon}
\end{figure}

Since the surface arising from the $12$-gon is too hard to visualize we display the mirror complex of the pentagon in
the Schlegel diagram of $\ncp{4}{5}$ in Figure~\ref{fig:M_5gon}.

\section{Concluding Remarks and Acknowledgments}

It can be shown that if the BBC sequence $\cT$ is polytopal then the cubical sphere $\bbc{\cT}$ is also
polytopal~\cite{RamanThiloGuenter}.  In view of Theorem~\ref{thm:polytopal} this gives the complete picture of the
construction as far as questions of polytopality are concerned.

For the visualization of the simplicial balls and the polyhedral surface $\mathcal{M}_{4,5}$ we used the software
packages \texttt{polymake}~\cite{MR1785292} and \texttt{JavaView}~\cite{javaview}.

We are indebted to Frank H. Lutz for his help in recognizing the sphere $\partial A_9$ in
Table~\ref{tab:BBCSequenceOfAlthulerSphere}.

\bibliographystyle{amsplain}
\bibliography{main}

\end{document}